\documentclass[12pt,reqno]{amsart}
\usepackage{xy, amssymb, amsfonts, amsbsy, amsthm, amsmath, amscd, latexsym, stmaryrd, epic, eepic,eucal}

\DeclareMathAlphabet{\mathpzc}{OT1}{pzc}{m}{it}

\newenvironment{dem}{\begin{proof}[\bf Proof]}{\end{proof}}

\newtheorem{theorem}{\bf Theorem}[section]
\newtheorem{lemma}[theorem]{\bf Lemma}
\newtheorem{propos}[theorem]{\bf Proposition}
\newtheorem{corol}[theorem]{\bf Corollary}
\newtheorem{claim}[theorem]{\bf Claim}

\theoremstyle{definition}
\newtheorem{defi}[theorem]{\bf Definition}
\newtheorem{oss}[theorem]{\bf Remark}

\newtheorem{exm}[theorem]{\bf Example}

\newcommand{\B}{\text{B}}

\newcommand{\C}{\mathbb C}

\newcommand{\F}{\mathcal F}
\newcommand{\Gm}{\mathbb G_{\textbf{m}}}
\newcommand{\Ga}{\mathbb G_{\textbf{a}}}

\newcommand{\M}{\mathfrak M}

\newcommand{\Ms}{\mathcal M}
\newcommand{\N}{\mathbb N}
\newcommand{\Nf}{\mathcal N}
\newcommand{\Of}{\mathcal O}
\newcommand{\Pro}{\mathbb P}
\newcommand{\Q}{\mathbb Q}

\newcommand{\Z}{\mathbb Z}
\newcommand{\aut}{\text{Aut}}

\newcommand{\az}{\gamma}

\newcommand{\cart}{\ar @{} [dr] |{\Box}}
\newcommand{\cat}{\mathpzc{Sch}_{\mathbb C}}

\newcommand{\cm}{\mathpzc c}
\newcommand{\cu}{\mathpzc C}

\newcommand{\hol}{\frac{\partial}{\partial z}}
\newcommand{\id}{\text{id}}

\newcommand{\km}{\mathpzc k}

\newcommand{\nm}{\mathpzc n}

\newcommand{\rat}{\otimes \mathbb Q}

\newcommand{\sla}{\mathfrak{sl}}
\newcommand{\spc}{\text{Spec} \mathbb C}
\newcommand{\spe}{\text{Spec}}
\newcommand{\spk}{\text{Spec}\Omega}

\newcommand{\unoa}{\xymatrix { *=0{\bullet} \ar@{-}[r] & *=0{\bullet} }}
\newcommand{\duea}{\xymatrix { *=0{\bullet} \ar@{-}[r] & *=0{\bullet} \ar@{-}[r] &*=0{\bullet} }}

\newcommand{\trea}{\xymatrix { *=0{\bullet} \ar@{-}[r] & *=0{\bullet} \ar@{-}[r] &*=0{\bullet} \ar@{-}[r] &*=0{\bullet} }}

\newcommand{\treb} {\xymatrix @R=8pt { & *=0{\bullet} \ar@{-}[dd] & \\ & &  \\ & *=0{\bullet} \ar@{-}[dl] \ar@{-}[dr] & \\ *=0{\bullet} & & *=0{\bullet}}}

\newcommand{\trebg}{\xymatrix @R=8pt { & *=0{} \ar@{-}[4,0] & \\ *=0{} \ar@{-}[0,2] & & *=0{} \\ *=0{} \ar@{-}[0,2] & & *=0{}\\ *=0{} \ar@{-}[0,2] & & *=0{}\\ & *=0{} &}}

\newcommand{\quattroa}{\xymatrix { *=0{\bullet} \ar@{-}[r] & *=0{\bullet} \ar@{-}[r] &*=0{\bullet} \ar@{-}[r] &*=0{\bullet} \ar@{-}[r] &*=0{\bullet} }}

\newcommand{\quattrob}{\xymatrix {& *=0{\bullet}  \ar@{-}[d]& \\ *=0{\bullet} \ar@{-}[r] & *=0{\bullet} \ar@{-}[r] \ar@{-}[d] & *=0{\bullet}\\ & *=0{\bullet} &}}

\newcommand{\quattroc}{\xymatrix @R=4pt {*=0{\bullet} \ar@{-}[dr] & & &\\ & *=0{\bullet} \ar@{-}[r] &*=0{\bullet} \ar@{-}[r] &*=0{\bullet}\\ *=0{\bullet} \ar@{-}[ur] & & &}}

\newcommand{\quattrocg}{\xymatrix @R=8pt { & *=0{} \ar@{-}[4,0] & & &\\ *=0{} \ar@{-}[0,2] & & *=0{} &*=0{} \ar@{-}[2,0] & \\ *=0{} \ar@{-}[0,4] & & & &*=0{}\\ *=0{} \ar@{-}[0,2] & & *=0{} &*=0{} &\\ & *=0{} & & &}}

\newcommand{\cinquenz}{\xymatrix { *=0{\bullet} \ar@{-}[d] & *=0{\bullet} \ar@{-}[d] \\ *=0{\bullet} \ar@{-}[r] & *=0{\bullet} \\
*=0{\bullet} \ar@{-}[u] & *=0{\bullet} \ar@{-}[u]}}

\input xy
\xyoption{all}

\begin{document}

\title[Tautological Classes of Rational Curves]{TAUTOLOGICAL CLASSES OF THE STACK OF RATIONAL NODAL CURVES}
\author{Damiano Fulghesu}
\address{Department of Mathematics,
University of Missouri,
Columbia, MO 65211}
\email{damiano@math.missouri.edu} 

\begin{abstract}
This is the second in a series of three papers in which we investigate the rational Chow ring of the stack $\M_{0}$ consisting of nodal curves of genus $0$. Here we define the basic classes: the classes of strata and the Mumford classes.
\end{abstract}

\maketitle

\section{Introduction}

\medskip

In the first paper \cite{fulg2} we described the stack $\M_0$, its stratification by nodes and more precisely by topological types given by dual graphs. When the graph $\Gamma$ corresponding to a curve $C$ has maximal multiplicity 3, we have an equivalence:
$$
\M_0^{\Gamma} \equiv \M_0 ^{\Gamma} \simeq \B \aut(C).
$$
This explicit description allows us to describe completely two kinds of classes.

{Strata classes:} for each tree $\Gamma$ with $\delta$ edges and maximal multiplicity 3, we get a class $\gamma_{\Gamma}$ of codimension $\delta$ in $\M_{0}$: the class of the closure of $\M_{0}^{\Gamma}$ in $\M_{0}$. We also compute the restriction of each $\gamma_{\Gamma}$ to all the rings $A^*(\M_{0}^{\Gamma'})\rat$ for any other tree $\Gamma'$ with maximal multiplicity 3 (Proposition \ref{refrestr}).

{Mumford classes:} as the map from the universal curve is not locally projective we cannot define these classes as in the usual sense because in Kresch's theory one has only pushforwards along projective morphisms.

However we are able to extend these classes in the following way (Section \ref{mumcla}): let $\cu \xrightarrow{\Pi} \M_{0}$ be the universal curve. Call $K \in A^1(\cu)\rat$ the first Chern class of the relative dualizing sheaf $\omega_{\cu/ \M_{0}}$, and set
   $$
   \km_{i} = \Pi_{*}(K^{i+1}).
   $$
Here we encounter a hard technical problem: the morphism $\cu \to \M_{0}$ is not projective (not even represented by schemes, as we have seen in \cite{fulg2}). 

If $n \leq 3$ then the pushforward $\Pi_{*}\omega^{\vee}_{\cu / \M_{0}}$ is a locally free sheaf of rank $3$ on the the open substack $\M_{0}^{\leq 3}$ (this fails for curves with $4$ nodes). Hence we have Chern classes  $\cm_{1}$, $\cm_{2}$ and $\cm_{3}$ of $\Pi_{*}\omega^{\vee}_{\cu / \M_{0}}$ in $A^*(\M_{0}^{\leq 3})$.

We have a pushforward $\Pi_{*} A^*(\mathcal C^{\Gamma})\rat \to A^*(\M_{0}^{\Gamma})\rat$ along the restriction $\mathcal C^{\Gamma} \xrightarrow{\Pi} \M_{0}^{\Gamma}$ of the universal curve, because the stacks involved are quotient stacks, and arbitrary proper pushforwards exists in the theory of Edidin and Graham. This allows, with Grothendieck-Riemann-Roch, to compute the restriction of the Mumford classes to each $A^*(\M_{0}^{\Gamma})\rat$ (Proposition \ref{defmumstrat}), even without knowing that the Mumford classes exist. It turns out that for each tree $\Gamma$ with at most $3$ nodes, the Mumford classes in $A^*(\M_{0}^{\Gamma})\rat$ are polynomials in the restrictions of $\cm_{1}$, $\cm_{2}$ and $\cm_{3}$, thought of as elementary symmetric polynomials in three variables. Then we define the classes $\km_{i} \in A^*(\M_{0}^{\leq 3})\rat$ as the suitable polynomials in the $\cm_{1}$, $\cm_{2}$ and $\cm_{3}$ (Definition \ref{defmumpol}).

Since $A^*(\M_{0}^{\leq3})\rat$ injects into the product of $A^*(\M_{0}^{\Gamma})\rat$ over trees with at most three nodes (Proposition \ref{injleqtre}), this gives the right definition.

\medskip

\section{Chow rings of strata}\label{stratclass}

In this paper we work on a fixed tree $\Gamma$ with maximal multiplicity 3 and $\delta$ edges. We indicate with $\Delta_1, \Delta_2, \Delta_3$ the sets of vertices that belongs respectively to one, two and three edges (and with $\delta_1, \delta_2, \delta_3$ their cardinalities).

We restrict the universal curve $\cu \to \M_0$ to $\M^{\Gamma}_0$
$$
\cu^{\Gamma} \xrightarrow{\Pi} \M^{\Gamma}_0.
$$
and consider the normalization (see \cite{Vis} Definition 1.18)
$$
\widehat{C}^{\Gamma} \xrightarrow{N} \cu^{\Gamma}.
$$

\begin{oss}
Given a curve $C \xrightarrow{\pi} T$ in $\M^{\Gamma}_0$ (that is to say a morphism $T \xrightarrow{f} \M^{\Gamma}_0$) we define $\widehat{C}$ as $T \times_{\M^{\Gamma}_0} \widehat{\cu}^{\Gamma}$. We have the following cartesian diagram
\begin{equation*}
\xymatrix{
\widehat{C} \cart \ar[r] \ar[d]_{n} & \widehat{\cu}^{\Gamma} \ar[d]^{N}\\
C \cart \ar[d]_{\pi} \ar[r] & \cu^{\Gamma} \ar[d]^{\Pi}\\
T \ar[r]^{f} & \M^{\Gamma}_0.
}
\end{equation*}

We notice that when $T$ is a reduced and irreducible scheme $\widehat{C} \xrightarrow{n} C$ is the normalization.

The map $\pi n: \widehat{C} \to T$ is proper as $\Pi N: \cu^{\Gamma} \to \M^{\Gamma}_0$ is. Then there exists a finite covering $\widetilde T \xrightarrow{} T$ (see \cite{Kn} Chapter 5, Theorem 4.1) such that we have the following commutative diagram
\begin{equation*}
\xymatrix{
 & \widehat C \ar[ddl]_g \ar[d]^n\\
 & C \ar[d]^\pi\\
\widetilde{T} \ar[r] &T
}
\end{equation*}
where the map $g$ has connected fibers.
\end{oss}

\begin{defi}
We define $\widetilde{\M}_0 ^{\Gamma}$ as the fibered category on $\cat$ whose objects are rational nodal curves $C \xrightarrow{\pi} T$ in $\M_0^{\Gamma}$ equipped with an isomorphism $\varphi: \coprod^{\Gamma}T  \to \widetilde T$ over $T$ and maps are morphisms in $\M_0^{\Gamma}$ that preserve isomorphisms.
\end{defi}

Straightforward arguments show the following

\begin{lemma}
The forgetful morphism
$$
\widetilde{\M}_0 ^{\Gamma} \to \M_0^{\Gamma}
$$
is representable finite \'etale and surjective.
\end{lemma}

We have a more explicit description of $\widetilde{\M}_0 ^{\Gamma}$: if we call $\Ms_{0,i}^n$ the stack of rational curves with $n$ nodes and $i$ sections, we can exhibit an equivalence
\begin{equation} \label{equivstrat}
\widetilde{\M}_0^{\Gamma} \cong \prod_{\alpha \in \Gamma} (\Ms_{0,e(\alpha)}^0)
\end{equation}
where, for each $\alpha \in \Gamma$, $e(\alpha)$. The proof of this is straightforward, however in the following we assume that the maximal multiplicity $k$ of $\Gamma$ is at most 3. In this case, as we have seen in  \cite{fulg2}:
$$
\M_0^{\Gamma} = \B \aut(C)
$$
where $\aut(C)$ is the group $\aut(\Gamma) \ltimes \left ( E^{\Delta_1} \times \Gm^{\Delta_2} \right )$. Let us call $\sigma$ the order of $\aut{\Gamma}$.

The \'etale covering of degree $\sigma$
\begin{eqnarray*}
\widetilde{\M}_0^{\Gamma} &\xrightarrow{\phi}& \M_0^{\Gamma}
\end{eqnarray*}
becomes
\begin{eqnarray*}
\B H &\xrightarrow{\phi}& \B \aut(C)
\end{eqnarray*}
where $H$ is the group $E^{\Delta_1} \times \Gm^{\Delta_2}$.

In order to compute the Chow ring on each stratum we only need equivariant intersection theory.

Let $\Gamma$ be a tree of maximal multiplicity at most 3 and let $C$ be the unique isomorphism class of curves of topological type $\Gamma$. Let us fix coordinates $[X,Y]$ on each component of $C$ such that
\begin{itemize}
\item on components with one node the point $[1,0]$ is the node,
\item on components with two nodes the points $[0,1]$ and $[1,0]$ are the nodes,
\item on components with three nodes the points $[0,1]$, $[1,1]$ and $[1,0]$ are the nodes.
\end{itemize}
Let us define on each component
\begin{eqnarray*}
0:=[0,1] \quad 1:=[1,1] \quad \infty:=[1,0]
\end{eqnarray*}

\begin{propos} \label{ringstratum}
Given $\Gamma$ as above, the Chow ring $A^*(\M^{\Gamma}_0)$ is
$$
A^*_{\aut(C)} \rat \cong (\Q[t_{\Delta_1}, r_{\Delta_2}])^{\aut(\Gamma)}
$$
where the action of an element $g \in \aut(\Gamma)$ on $\Q[x_{\Delta_1}, y_{\Delta_2}]$ is the obvious permutation on the $\Delta_1 \cup \Delta_2$ variables together with multiplication by (-1) of $r$-variables corresponding to components of which $g$ exchanges $0$ and $\infty$.
\end{propos}

\begin{dem}
The group $E^{\Delta_1} \times \Gm^{\Delta_2}$ is a normal subgroup of $\aut(C)$.

So we can apply the following

\begin{lemma} \label{isophi} \cite{Ve}
Given an exact sequence of algebraic groups over $\C$
$$
\xymatrix{
1 \ar[r] & H  \ar[r]^\phi & G \ar[r]^\psi &  F \ar[r] & 1.
}
$$
with $F$ finite and $H$ normal in $G$, we have
$$
A^*_{G} \rat \cong (A^*_{H} \rat )^{F}
$$
\end{lemma}
and conclude that
$$
A^*_{\text{Aut}(C)} \rat \cong (A^*_{ E^{\Delta_1}\times (\Gm)^{\Delta_2}} \rat)^{\text{Aut}(\Gamma)}
$$

So we have reduced to compute $A^*_{E^{\Delta_1} \times (\Gm)^{\Delta_2}}$ and the action of $\aut(\Gamma)$ on it.

\begin{claim} \label{pr}
The ring $A^*_{E^{\Delta_1} \times (\Gm)^{\Delta_2} } \rat$ is $\mathbb Q[x_{\Delta_1}, y_{\Delta_2}]$, that is to say it is algebraically generated by  $\Delta_1 \cup \Delta_2$ independent generators of degree 1.
\end{claim}

We use the following fact (\cite{Ve} Proposition 2.8):
\begin{lemma}
For every linear algebraic group $G$, we have
$$
A^*_{G \times \Gm} \cong A^*_{G} \otimes_{\Z} A^*_{\Gm}.
$$
\end{lemma}
By recalling that $ A^*_{\Gm} \rat$ is isomorphic to $\Q [r]$, where $r$ is an order one class, we have
$$
A^*_{E^{\Delta_1} \times (\Gm)^{\Delta_2} } \rat \cong A^*_{E^{\Delta_1}} \otimes_{\Q} \Q[y_{\Delta_2}].
$$
Using the fact that the group $E$ is the semidirect product
$$
\xymatrix{
0 \ar[r] & \Ga  \ar[r]^{\varphi} & E \ar@<-1ex>[r]_{\rho} & \Gm \ar@<-1ex>[l]_{\psi} \ar[r] & 1,
}
$$
we can explicit that
$$
A^*_{E^{\Delta_1}} \rat \cong \Q[x_{\Delta_1}].
$$

Now we describe the action of $\aut(\Gamma)$ on
$$
A^*_{E^{\Delta_1} \times \Gm^{\Delta_2}} \cong \Q[x_{\Delta_1}, y_{\Delta_2}].
$$
In order to do this we need a more explicit description of the classes $x_{\Delta_1}, y_{\Delta_2}$ through the equivalence (\ref{equivstrat})

\begin{eqnarray*}
(\Ms^0_{0,1})^{\Delta_1} \times (\Ms^0_{0,2})^{\Delta_2} \times (\Ms^0_{0,3})^{\Delta_3} &\cong& \widetilde{\M}_0^{\Gamma}.
\end{eqnarray*}
Since $\Ms^0_{0,3} \simeq \spe \C$ we have
\begin{eqnarray*}
(\Ms^0_{0,1})^{\Delta_1} \times (\Ms^0_{0,2})^{\Delta_2} &\xrightarrow{\phi}& \M_0^{\Gamma}.
\end{eqnarray*}
Moreover we have
\begin{eqnarray*}
\Ms^0_{0,1} &\simeq& \B E\\
\Ms^0_{0,2} &\simeq& \B \Gm.
\end{eqnarray*}

Let $\alpha$ be a vertex of $\Gamma$ such that $e(\alpha)=1$ or $2$. On the component $\Ms_{0,e(\alpha)}^0$ let us consider the universal curve
$$
\cu^{\alpha} \xrightarrow{\widetilde{\Pi}} \Ms_{0,e(\alpha)}^0. 
$$
On $\Pro^1_{\C}$ we fix coordinates and we define the points
\begin{eqnarray*}
z_\infty &:=& [1,0]\\\
z_0 &:=& [0,1]
\end{eqnarray*}
We can write
\begin{itemize}
\item $\cu^{\alpha} \simeq [\Pro^1_{\C}/E]$ when $e(\alpha)=1$
\item $\cu^{\alpha} \simeq [\Pro^1_{\C}/\Gm]$ when $e(\alpha)=2$
\end{itemize}
Let us consider on $\Pro^1_{\C}$ the linear bundles $\Of(z_\infty)$ and $\Of(z_0)$. We have a natural action of $\Gm$ on global sections of both of them induced by the action of $\Gm$ (that, we recall, fixes $z_\infty$ and $z_0$) on $\Pro^1_{\C}$. Similarly we have an action of the group $E$ on global sections of $\Of(z_\infty)$.
When $e(\alpha)=1$ set
\begin{equation*}
\psi^{1}_{\infty,\alpha}:=c^E_1(H^0(\Of(z_\infty), \Pro^1_{\C}));
\end{equation*}
while, if $e(\alpha)=2$ set
\begin{eqnarray*}
\psi^{2}_{\infty,\alpha} &:=& c^{\Gm}_1(H^0(\Of(z_\infty), \Pro^1_{\C}))\\
\psi^{2}_{0,\alpha} &:=& c^{\Gm}_1(H^0(\Of(z_0), \Pro^1_{\C})).
\end{eqnarray*}
Clearly, when $e(\alpha)=2$, we have
$$
\psi^2_{\infty,\alpha} + \psi^2_{0,\alpha}=0
$$
Now we define
\begin{eqnarray*}
&& t_{\alpha}:=\psi^1_{\infty,\alpha} \text{ when $e(\alpha)=1$}\\
&& r_{\alpha}:=\psi^2_{\infty,\alpha} \text{ when $e(\alpha)=2$}
\end{eqnarray*}
We can write for each $\alpha$ such that $e(\alpha)=2$
$$
r_\alpha=\frac{\psi^2_{\infty,\alpha} - \psi^2_{0,\alpha} }{2}.
$$
Clearly all the classes $t_{\Delta_1}$ and $r_{\Delta_2}$ are of order one and independent. From what we have seen above these classes generates the ring $A^*(\widetilde{\M}^{\Gamma}_0) \rat$ and we have
$$
A^*(\widetilde{\M}^{\Gamma}_0) \rat \simeq \Q[t_{\Delta_1}, r_{\Delta_2}].
$$
Now we can describe the action of $\aut{\Gamma}$ on $\Q[t_{\Delta_1}, r_{\Delta_2}]$. An element $g \in \aut(\Gamma)$ acts on $C$ with a permutation $g_1$ on the components with one node and a permutation $g_2$ of components with two nodes. As we have chosen coordinates on $C_0$ such that $\infty$ corresponds to the node of the terminal components, we make $g_1$ act directly to the set $\{ t_{\Delta_1} \}$. We make $g_2$ act similarly on the set $\{ r_{\Delta_2} \}$ but we have in addition to consider the sign, that is to say that when $g_2$ sends a vertex $P$ of $\Gamma$ to another vertex $\beta$ (such that $e(\alpha)=e(\beta)=2$), we have two possibilities
\begin{itemize}
\item the automorphism $g$ exchange coordinates 0 and $\infty$ and so we have
$$
g(r_{\alpha})=g \left( \frac{\psi^2_{\infty,\alpha} - \psi^2_{0,\alpha}}{2}\right ) = \frac{\psi^2_{0,\beta} - \psi^2_{\infty,\beta}}{2} = -r_{\beta}
$$
\item the automorphism $g$ sends 0 in 0 and $\infty$ in $\infty$; in this case we have
$$
g(r_{\alpha}) = r_\beta.
$$
 \end{itemize}

\end{dem}

\begin{defi}
We define  $\gamma_i$ as the class in  $A^*(\M_0^{\leq i})$ of $\M_0^i$. We will indicate with $\gamma_i \in A^*(\M_0)$ also the class of the closure of $\M_0^i$ in $\M_0$. Similarly we define $\gamma_{\Gamma}$ as the class of the closure of $\M^{\Gamma}_0$ in $\M_0$.
\end{defi}

\begin{propos}\label{classstratum}
Let $\Gamma$ be a tree of maximal multiplicity at most three (except the single point) and $C$ be the curve of topological type $\Gamma$. Let us consider the \'etale covering
$$
\widetilde{\M}_0^\Gamma \xrightarrow{\phi} \M_0^{\Gamma}
$$
Let $C_{\Gamma}$ be the components of $C$ (which we see as vertex of $\Gamma$). Let $E(\Gamma)$ be the set of edges. If $(\alpha, \beta) \in E(\Gamma)$ we call $z_{\alpha \beta}$ the common point of $C_\alpha$ and $C_\beta$.
Then, by following notation of Proposition \ref{ringstratum}, we have
$$
\phi^* \gamma_{\Gamma} = \prod_{(\alpha, \beta) \in E(\Gamma)} \left( \psi^{e(\alpha)}_{z_{\alpha \beta}|_{\alpha}, \alpha} + \psi^{e(\alpha)}_{z_{\alpha \beta}|_{\beta}, \beta} \right).
$$
In particular the classes of each stratum of $\M^{\leq 3}_0$
after fixing coordinates on $C$ and ordering components of $\Delta_1$ and $\Delta_2$, are:
\medskip

\begin{center}
\begin{tabular}{|c|c|}
\hline
\small{Graph $(\Gamma)$} & class of stratum $\phi^*\gamma_{\Gamma}$\\
\hline
\small{$$ \unoa $$} & $t_1 + t_2$\\
\hline
\small{$$ \duea $$} & $(t_1 - r_1)(t_2 + r_1)$ \\
\hline
\small{$$ \trea $$} & $(t_1 - r_1)(r_1 + r_2)(t_2 - r_2)$\\
\hline
\small{$$ \treb $$} & $t_1 t_2 t_3$\\
\hline
\end{tabular}
\end{center}

\end{propos}

\begin{dem}
For every tree $\Gamma$ of maximal multiplicity at most three (except the single point), let us consider the regular embedding
$$
\M_0^\Gamma= \B \aut(C) \xrightarrow{in} \M_0^{\leq \delta}
$$
where $\delta$ is the number of edges of $\Gamma$.

Let us consider the normal bundle
$$
N_\Gamma:= N_{\M_0^\Gamma / \M_0^{\leq \delta}}.
$$
It is known that $N_\Gamma=\text{def}_{\Gamma}$ is the space of first order deformations of $\M_0^{\Gamma}$
$$
{\bigoplus_{(\alpha,  \beta) \in E(\Gamma)}} T_{z_{\alpha \beta}}(C_\alpha) \otimes T_{z_{\alpha \beta}}(C_\beta).
$$

As usual let us consider the \'etale covering
\begin{eqnarray*}
&& \widetilde \M_0^{\Gamma} \xrightarrow{\phi} \M_0^{\Gamma}\\
&& \B H \to \B ( \aut (\Gamma) \ltimes H )
\end{eqnarray*}
and set $\widetilde N _{\Gamma}= \phi^* N_{\Gamma}$. We notice that the point $z_{\alpha\beta}$ on each component $C_\alpha$ (which we call $z_{\alpha\beta}|_{\alpha}$) is 0, 1 or $\infty$.

By using notation of Proposition \ref{ringstratum} we have on $\widetilde \M^{\Gamma}_0$
$$
c^{G_{\alpha}}_1(T_{z_{\alpha \beta}}(C_\alpha))= \psi^{e(\alpha)}_{z_{\alpha \beta}|_{\alpha}, \alpha}
$$
where $G_{\alpha}=E$ if $e(\alpha)=1$, $G_{\alpha}=\Gm$ if $e(\alpha)=2$ and $G_{\alpha}=\id$ if $e(\alpha)=3$
Notice that $c^{G_{\alpha}}_1(T_{z_{\alpha \beta}}$ is zero when $e(\alpha)=3$, consequently
$$
c^{H}_{\text{top}}(\widetilde{N}_{\Gamma})= \prod_{(\alpha, \beta) \in E(\Gamma)} \left( \psi^{e(\alpha)}_{z_{\alpha \beta}|_{\alpha}, \alpha} + \psi^{e(\beta)}_{z_{\alpha \beta}|_{\beta}, \beta} \right).
$$
We have the following relation in $A^*(\M^{\Gamma}_0)$

\begin{equation} \label{nor}
in^*in_*[\M_0^{\Gamma}]=c^{H}_{top}(N_\Gamma) \cap [\M_0^{\Gamma}].
\end{equation}
and so
$$
\phi^*\gamma_{\Gamma}= c^{H}_{\text{top}}(\widetilde{N}_{\Gamma})
$$
\end{dem}

\begin{oss} \label{zerodiv}

Given a choice of coordinates the class $\phi^* \gamma_{\Gamma}$ is invariant for the action of $\aut (\Gamma)$ given in Proposition \ref{ringstratum}, so we actually can see it as the class $\gamma_{\Gamma}$ in $\M^\Gamma_0$.

Moreover we have shown that these classes are not 0-divisor in the ring $A^*(\M^{\Gamma}_0)$ for each $\Gamma$ corresponding to a stratum in $\M^{\leq 3}_0$.

This fails if we consider more than four nodes. For example if we consider the following graph $\Gamma$
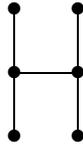
\begin{equation*}
\cinquenz
\end{equation*}
we have $\phi^* \gamma_{\Gamma}=0$.
\end{oss}

For future reference we can now state the following

\begin{propos} \label{injleqtre}
The Chow ring $A^*(\M^{\leq 3}_0) \rat$ injects into the product of $A^*(\M^{\Gamma}_0)\rat$ over trees with at most three edges.
\end{propos}
\begin{dem}
From Proposition \ref{classstratum} and Remark (\ref{zerodiv}), we have, for each $\delta \leq 3$, the following exact sequence of additive groups
$$
0 \to A^*(\M^{\delta}_0) \rat \xrightarrow{i^{\delta}_*} A^*(\M^{\leq \delta}_0) \rat \xrightarrow{j^{*\delta}} A^*(\M^{\leq (\delta -1)}_0)
$$
where
\begin{eqnarray*}
&&i^{\delta}: \M^{\delta}_0 \to \M^{\leq \delta}_0\\
&&j^{\delta}: \M^{\leq (\delta -1)}_0 \to \M^{\leq \delta}_0
\end{eqnarray*}
are the natural closed embeddings.

Let us consider the morphism
$$
A^*(\M^{\leq 3}_0) \rat \xrightarrow{\psi} \prod^{3}_{\delta=0} A^*(\M^{\delta}_0) \rat
$$
as the product of the maps $i^{3*}$, $i^{2*} j^{3*}$, $i^{1*} j^{2*}j^{3*}$ and $j^{1*}j^{2*}j^{3*}$.
Let $a$ be an element of $A^*(\M^{\leq 3}_0) \rat$ different from zero. If $\psi(a)$ is zero then it cannot be in the image of $i^3_*$ consequently $j^{3*}(a) \in A^*(\M^{\leq 2}_0)$ is different from zero. We can continue till we obtain that $j^{1*}j^{2*}j^{3*}$ is different from zero: absurd.
\end{dem}

\section{Restriction of classes to strata} \label{restrcla}
\medskip

Let us consider two trees $\Gamma$ and $\Gamma'$ with maximal multiplicity at most $3$ and number of edges respectively equal to $\delta$ and $\delta'$.

\begin{defi}
Given two graphs as above we call an ordered deformation of $\Gamma$ into $\Gamma'$ any surjective map of vertices $d: \Gamma' \to \Gamma$ such that
\begin{enumerate}
\item for each $P,Q \in \Gamma'$ we have $d(P) = d(Q)=A \in \Gamma'$ only if for each $R$ in the connected path from $P$ to $Q$ we have $d(R)=A$;
\item for each edge $(P,Q) \in \Gamma'$ such that $d(P) \neq d(Q)$ there must be an edge in $\Gamma$ between $d(P)$ and $d(Q)$.
\end{enumerate}
We denote by $\text{def}_o(\Gamma, \Gamma')$ the set of deformations.
\end{defi}

\begin{exm} \label{exmdef}
Let $\Gamma$ and $\Gamma'$ be the following graphs
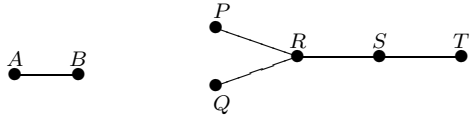
\begin{equation*}
\vcenter{\xymatrix { *=0{\bullet} \ar@{-}[r]^<A & *=0{\bullet} \ar@{}[r]^< B& *=0{}  }}
\qquad
\vcenter{\xymatrix @R=4pt {*=0{\bullet} \ar@{-}[dr]^<P & & & &\\
& *=0{\bullet} \ar@{-}[r]^<R &*=0{\bullet} \ar@{-}[r]^<S &*=0{\bullet} \ar@{}[r]^< T & *=0{}\\
*=0{\bullet} \ar@{-}[ur]_<Q & & & &}}
\end{equation*}

we have the following 8 ordered deformations

\begin{small}
\begin{tabular}{|ll|ll|}
\hline
$1) (P,R,S,T) \mapsto A$ & $Q \mapsto B$ & $5) Q  \mapsto A$ & $(P,R,S,T)  \mapsto B$\\
\hline
$2) (Q,R,S,T) \mapsto A$ & $P \mapsto B$ & $6) P  \mapsto A$ & $(Q, R, S, T)  \mapsto B$\\
\hline
$3) (P,Q,R) \mapsto A$ & $(S,T) \mapsto B$ & $7) (S,T)  \mapsto A$ & $(P, Q, R)  \mapsto B$\\
\hline
$4) (P,Q,R,S) \mapsto A$ & $T \mapsto B$ & $8) T  \mapsto A$ & $(P, Q, R, S) \mapsto B$\\
\hline
\end{tabular}
\end{small}
\end{exm}

There exist two different equivalence relations in $\text{def}_o(\Gamma, \Gamma')$.
We say that two elements $d_1,d_2$ are in $\sim_\Gamma$ if there exists a $\gamma \in \aut(\Gamma)$ such that $d_2= \gamma d_1$. Similarly we say that two elements $d_1,d_2$ are in $\sim_{\Gamma'}$ if there exists a $\gamma' \in \aut(\Gamma')$ such that $d_2= d_1\gamma'$.

\begin{defi}
We call $\Gamma-$deformations (or simply deformations) from $\Gamma$ to $\Gamma'$ the set
$$
\text{def}_{\Gamma}(\Gamma, \Gamma') := \text{def}_o(\Gamma, \Gamma')/ \sim_{\Gamma}.
$$
We call $\Gamma'-$deformations from $\Gamma$ to $\Gamma'$ the set
$$
\text{def}_{\Gamma'}(\Gamma, \Gamma') := \text{def}_o(\Gamma, \Gamma')/ \sim_{\Gamma'}.
$$

\end{defi}

In the above example we take as representatives of $\Gamma-$deformations the first 4 ordered deformations. On the other hand we have $$1 \sim_{\Gamma'} 2 \qquad 5 \sim_{\Gamma'} 6$$  

From topological arguments we have the following
\begin{propos}\label{defclaim}
Let $C \xrightarrow{\pi} T$ be a family of rational nodal curves over an irreducible scheme $T$. Suppose further that the generic fiber has  topological type $\Gamma$, then there exists a fiber of topological type $\widehat{\Gamma}$ only if there exists an ordered deformation of $\Gamma$ into $\widehat{\Gamma}$.
\end{propos}

Now let us consider the \'etale map
$$
\widetilde{\M}^{\Gamma}_0 \xrightarrow{\phi} \M^{\Gamma}_0.
$$
We have given above a description of $A^*{\M_0}$ as the subring of polynomials of $A^*(\widetilde{\M}^{\Gamma}_0)$ in the classes (corresponding to sections of $\widetilde{\M}_0^{\Gamma}$) $t_1, \dots, t_{\delta_1}, r_1, \dots, r_{\delta_2}$ invariant for the action of $\aut(\Gamma)$.

We call $\Ms_{0,i}$ the stack of rational nodal curves with $i$ sections.
Let $\left( \widetilde{\M}^{\Gamma}_0 \right)^{\leq \delta' -\delta}$ be the substack of
$$
\left( \Ms_{0,1} \right)^{\Delta_1} \times \left( \Ms_{0,2} \right)^{\Delta_2} \times \left( \Ms_{0,3}\right)^{\Delta_3}
$$
whose fibers have at most $\delta' -\delta$ nodes (the sum of nodes is taken over all the connected components).
Polynomials in $\Q[t_1, \dots, t_{\delta_1}, r_1, \dots, r_{\delta_2}]$ has a natural extension to $\left( \widetilde{\M}^{\Gamma}_0 \right)^{\leq \delta' -\delta}$. Let us fix one of such polynomials $a$ which are invariants for the action of $\aut{\Gamma}$.

The \'etale covering
$$
\widetilde{\M}^{\Gamma}_0 \xrightarrow{\phi} \M^{\Gamma}_0
$$
is obtained by gluing sections in a way which depends on $\Gamma$.

By gluing sections in the same way we obtain a functor
$$
\left( \widetilde{\M}^{\Gamma}_0 \right)^{\leq \delta' -\delta} \xrightarrow{\Pi} \M^{\delta'}_0.
$$

\begin{corol}\label{cordefclaim}
With the above notation we have that the closure $\overline \M^{\Gamma}_{0}$ of $\M^{\Gamma}_0$ in $\M^{\leq \delta'}_0$ is
$$
\Pi \left( \left( \widetilde{\M}^{\Gamma}_0 \right)^{\leq \delta' -\delta} \right).
$$
\end{corol}
\begin{dem}
From Proposition \ref{defclaim} we have
$$
\overline \M^{\Gamma}_0 \subseteq \Pi \left( \left( \widetilde{\M}^{\Gamma}_0 \right)^{\leq \delta' -\delta} \right).
$$
On the other hand let $C_{\Omega} \xrightarrow{\pi} \spk$ be the image in $\M^{\delta'}_0$ of a geometric point of $\left( \widetilde{\M}^{\Gamma}_0 \right)^{\leq \delta' -\delta}$. The dual graph $\Gamma_{\Omega}$ of $C_{\Omega}$ is a deformation of $\Gamma$. In order to show that $C_{\Omega} \xrightarrow{\pi} \spk$ is a geometric point of $\overline{\M}_0^{\Gamma}$, we fix a deformation $d: \Gamma_{\Omega} \to \Gamma$. For each vertex $A$ of $\Gamma$, the set $d^{-1}(A)$ is a subtree of $\Gamma_{\Omega}$. We can give a deformation $C_A \xrightarrow{\pi} T$ of $C_{\Omega}$ such that the generic fiber is $\Pro^1$. Furthermore we can define on $C_A \xrightarrow{\pi} T$ a family of $E(A)$ that respect $d$. At last we glue all $C_A$ along sections and obtain a deformation $C \xrightarrow{\pi} T$ of  $C_{\Omega} \xrightarrow{\pi} \spk$ in $\overline{\M}_0^{\Gamma}(\spk)$.
\end{dem}

\begin{defi}
We call the image of $\Pi$:
$$
\M_0^{\text{def}(\Gamma, \delta')}:=  \Pi \left( \left( \widetilde{\M}^{\Gamma}_0 \right)^{\leq \delta' -\delta} \right) \subset \M_0^{\leq \delta'}.
$$
\end{defi}

\begin{propos}
The map
$$
\Pi: \left( \widetilde{\M}_0^{\Gamma}\right)^{\delta'- \delta} \to \M_0^{\leq \delta'}
$$
is finite.
\end{propos}
\begin{dem}
({\bf sketch}) We have to prove that $\Pi$ is representable, with finite fibers and proper. The not trivial property to verify is properness. We can prove it through the valutative criterion (see \cite{hrt} p.101).
\end{dem}
Therefore the map $\Pi$ is finite hence projective, so we have the push-forward \cite{kre}
$$
\Pi_*: A^* \left( \left( \widetilde{\M}_0^{\Gamma}\right)^{\delta'- \delta} \right) \to A^* \left( \M_0^{\text{def}(\Gamma, \delta')}\right)
$$

\begin{defi}
Let $\Gamma$ be a tree with maximal multiplicity $\leq 3$ and $\Gamma'$ a deformation of $\Gamma$ (with maximal multiplicity $\leq 3$) with $\delta'$ edges.
Let $a$ be a class in $A^*(\M^{\Gamma}_0) \rat$ and $\tilde a$ its $\delta'$-lifting.
With reference to the cartesian diagram
\begin{equation*}
\xymatrix@=3pc{
\M^{\Gamma'}_0 \times_{\M^{\delta'}_0} \left ( \widetilde{\M}^{\Gamma}_0 \right )^ {\delta' -\delta} \cart \ar[r]^-{pr_2} \ar[d]_{pr_1} &\left ( \widetilde{\M}_0^{\Gamma} \right ) ^{\leq \delta' - \delta} \ar[d]^{\Pi}\\
\M_0^{\Gamma'} \ar[r]^{in} &\M_0^{\leq \delta'}
}
\end{equation*}
we define
$$
\Psi(\Gamma, \Gamma'):= \M^{\Gamma'}_0 \times_{\M^{\delta'}_0} \left ( \widetilde{\M}^{\Gamma}_0 \right )^ {\delta' -\delta}.
$$
\end{defi}

Topological arguments show the following
\begin{propos}\label{diagstrat}
$\Psi(\Gamma, \Gamma')$ is a disjoint union of components that we write as
$$
\Psi(\Gamma, \Gamma')= :\coprod_{\xi \in \text{def}_{\Gamma'}(\Gamma, \Gamma')} \Psi(\Gamma, \Gamma')_\xi
$$
where the union is taken over the set of ordered deformations up to $\sim_{\Gamma'}$.
\end{propos}

\begin{exm} $\;$

Let $\Gamma$ and $\Gamma'$ be the graphs of the Example \ref{exmdef}. We have $\delta'-\delta=3$ and
$$
\widetilde{\M}^{\Gamma}_0= \Ms^0_{0,1} \times \Ms^0_{0,1}
$$
with a double covering
$$
\widetilde{\M}^{\Gamma}_0 \xrightarrow{\phi} \M^{\Gamma}_0.
$$
Clearly $\Psi(\Gamma, \Gamma')$ is an inclusion of components
in 
\begin{equation}
\coprod _{i,j: i+j=3} \left( \Ms^i_{0,1} \times \Ms^j_{0,1} \right)
\end{equation}

We have that $\Psi(\Gamma, \Gamma')$ has 6 components which corresponds to deformations (enumerated in Example \ref{exmdef}) up to $\sim_{\Gamma'}$, with the following inclusions:
\begin{itemize}
\item deformation 1 (which is $\Gamma'$-equivalent to 2) and 4 correspond to two connected components of $\Ms^{3}_{0,1} \times \Ms^{0}_{0,1}$
\item deformation 3 corresponds to a connected component of $\Ms^{2}_{0,1} \times \Ms^{1}_{0,1}$
\item deformation 7 corresponds to a connected components of $\Ms^{1}_{0,1} \times \Ms^{2}_{0,1}$
\item deformation 5 (which is $\Gamma'$-equivalent to 6) and 8 corresponds to two connected components of $\Ms^{0}_{0,1} \times \Ms^{3}_{0,1}$
\end{itemize}
\end{exm}

\begin{defi}
Let us consider a class $a$ in $A^*(\M^{\Gamma}_0) \rat$.
With reference to the \'etale covering
$$
\widetilde{\M}^{\Gamma}_0 \xrightarrow{\phi} \M^{\Gamma}_0
$$
we define the {\it lifting} of $a$
$$
\widetilde{a}:= \frac{\phi^*a}{\sigma}\in A^*(\widetilde{\M}^{\Gamma}_0) \rat,
$$
this means that $\phi_*(\widetilde{a})=a$.
Since $\widetilde a$ can be written as a polynomial in the classes $\psi$ defined in Proposition \ref{ringstratum}
that depends only on $\Gamma$, it has a natural extension to $A^*\left( \left( \widetilde{\M}_0^{\Gamma}\right)^{\delta'- \delta} \right) \rat$ that we call $\delta'$-lifting of $a$ and we still write $\tilde a$.
\end{defi}

For every $\xi \in \text{def}_{\Gamma'}(\Gamma, \Gamma')$ we have the related commutative diagram
\begin{equation*}
\xymatrix@=3pc{
\Psi(\Gamma, \Gamma')_\xi \ar[r]^{pr_2^\xi} \ar[d]_{pr_1^\xi} & \left ( \widetilde{\M}_0^{\Gamma} \right ) ^{\leq \delta' - \delta} \ar[d]^{\Pi}\\
\M_0^{\Gamma'} \ar[r]^{in} &\M_0^{\leq \delta'}
}
\end{equation*}
The map $pr_2^\xi$ is a closed immersion of codimension $\delta' - \delta$. Let us still call $\widetilde a$ the pullback of the polynomial $\widetilde a$ through $pr_2^{\xi}$.
By the excess intersection formula (Section 6.3 \cite{ful}, similar arguments show it for algebraic stacks)
we have
$$
in^!(\widetilde a)= \sum_{\xi \in \text{def}_{\Gamma'}(\Gamma, \Gamma')} (\widetilde a \cdot c_{top}[\Nf^{\xi}]).
$$
where $\Nf^{\xi}:=(pr_1^{\xi*} \Nf_{in})/ \Nf_{pr_2^{\xi}}$.

For each $\xi \in \text{def}_{\Gamma'}(\Gamma, \Gamma')$ we have an \'etale covering
$$
f_{\xi}: \widetilde \M^{\Gamma'} \to \Psi(\Gamma,\Gamma')_{\xi}
$$
that glue along sections.

For every $\xi \in \text{def}_{\Gamma'}(\Gamma, \Gamma')$ let us call $\sigma'_{\xi}$ the degree of $pr^{\xi}_1$.

We have the following commutative diagram
\begin{equation*}
\xymatrix{
\widetilde{\M}_0^{\Gamma'} \ar[rr]^{f_{\xi}} \ar[dr]_{\phi'} && \Psi(\Gamma, \Gamma')_\xi \ar[dl]^{pr_1^\xi}\\
& \M_0^{\Gamma'} &
}
\end{equation*}
from which we have
$$
\text{ord} \phi' = \text{ord} f_{\xi} \cdot \text{ord} \sigma'_{\xi}.
$$
\begin{oss} \label{restrpol}
The order of $f_{\xi}$ is the number of $g \in \aut(\Gamma')$ such that for each deformation $d$ associated to $\xi$ the deformation $d \circ g$ is $\Gamma'$-equivalent to $d$.
\end{oss}
We want to explicit ${\phi'}^*pr_{1*}in^!(\widetilde a)$ in $A^* ( \widetilde{\M}_0^{\Gamma'} ) \rat$: being invariant for the action of $\aut(\Gamma')$ we can see it as a class in $\M^{\Gamma'}_0$ which is the restriction of the extension of the class $a \in A^* \left( \M^{\Gamma}_0 \right) \rat$.

\begin{propos}\label{restriction}
We have
\begin{equation*}
{\phi'}^*pr_{1*}in^!(\widetilde a)=\sum_{\xi \in \text{def}_{\Gamma'}(\Gamma, \Gamma')} \sigma'_{\xi} f^*_{{\xi}}(\widetilde a) \cdot c_{top}(\widetilde{\Nf}^{\xi})
\end{equation*}
where
$$
\widetilde{\Nf}^{\xi}:= f^*_{{\xi}} {\Nf}^{\xi}.
$$
\end{propos}

\begin{dem}
Putting together all the above remarks and definitions we have
\begin{eqnarray*}
{\phi'}^*pr_{1*}in^!(\widetilde a)&=& {\phi'}^*pr_{1*} \sum_{\xi \in \text{def}_{\Gamma'}(\Gamma, \Gamma')} (\widetilde a \cdot c_{top}[\Nf^{\xi}])\\
&=&  {\phi'}^* \sum_{\xi \in \text{def}_{\Gamma'}(\Gamma, \Gamma')} pr^{\xi}_{1*}(\widetilde a \cdot c_{top}[\Nf^{\xi}])\\
&=& \sum_{\xi \in \text{def}_{\Gamma'}(\Gamma, \Gamma')} f^{*}_{\xi} pr^{\xi*}_{1} pr^{\xi}_{1*}(\widetilde a \cdot c_{top}[\Nf^{\xi}])\\
&=& \sum_{\xi \in \text{def}_{\Gamma'}(\Gamma, \Gamma')} \sigma'_{\xi} f^*_{{\xi}}(\widetilde a) \cdot c_{top}(\widetilde{\Nf}^{\xi})
\end{eqnarray*}
\end{dem}

We explicit now the computation of classes corresponding to strata.

\begin{propos}\label{refrestr}
Given a tree $\Gamma$ with maximal multiplicity $\leq 3$, let $\gamma_{\Gamma}$ be the class of $\overline \M^{\Gamma}_0$ in $\M_0$ and let $\Gamma'$ be another tree with maximal multiplicity $\leq 3$.

If $\Gamma'$ is a deformation of $\Gamma$, then the restriction of $\gamma_{\Gamma}$ to $A^* \left( \M^{\Gamma'} \right) \rat$ is (following the above notation)
\begin{equation*}
\sum_{\xi \in \text{def}_{\Gamma'}(\Gamma, \Gamma')} \frac{\sigma'_{\xi}}{\sigma}  c_{top}(\widetilde{\Nf}^{\xi}).
\end{equation*}

If $\Gamma'$ is not a deformation of $\Gamma$ then the restriction is 0.
\end{propos}

\begin{dem}
In case $\Gamma'$ is a deformation of $\Gamma$, the polynomial $a$ of Proposition (\ref{restriction}) is 1, consequently $\widetilde a$ is $1/\sigma$. With reference to the \'etale covering $\widetilde{\M}^{\Gamma'}_0 \xrightarrow{\phi} \M^{\Gamma'}_0$ and the inclusion $\M^{\Gamma'}_0 \xrightarrow{in} \M^{\delta'}_0$, from Proposition \ref{restriction} we obtain
\begin{equation*}
\sum_{\xi \in \text{def}_{\Gamma'}(\Gamma, \Gamma')} \frac{\sigma'_{\xi}}{\sigma}  c_{top}(\widetilde{\Nf}^{\xi}).
\end{equation*} 
If $\Gamma'$ is not a deformation of $\Gamma$ then using Corollary \ref{cordefclaim} and basic topological arguments, we have that $\M^{\Gamma'}_0$ does not intersect the closure of $\M^{\Gamma}_0$ in $\M_0$, so the restriction of the class must be 0.
\end{dem}

\begin{oss}
For each $\xi \in \text{def}_{\Gamma'}(\Gamma, \Gamma')$, we have an exact sequence of sheaves
$$
0 \to f^*_{\xi} \Nf_{pr^{\xi}_2} \to  {\phi'}^* \Nf_{in} \to \widetilde \Nf^{\xi} \to 0.
$$
We have seen in Proposition \ref{classstratum} how to compute $c_{top} ({\phi'}^* \Nf_{in})$, similarly we can compute $c_{top} (f^*_{\xi} \Nf_{pr^{\xi}_2})$ and finally we consider the following relation (that follows from the exact sequence)
$$
c_{top} ({\phi'}^* \Nf_{in})= c_{top} (f^*_{\xi} \Nf_{pr^{\xi}_2}) \cdot c_{top} (\widetilde \Nf^{\xi}).
$$
\end{oss}
We carry on the calculation for trees with at most three nodes in the last Section.

\section{Mumford classes}\label{mumcla}

Given a tree $\Gamma$ with at most four vertices, we have the restriction of the dualizing sheaf $\omega_0 := \omega_{\cu/ \M_0}$ on the universal curve $\cu^{\Gamma}$ of $\M^{\Gamma}_0$. We  call $\omega^{\Gamma}_0:=\omega_{\cu^{\Gamma} / \M^{\Gamma}_0}$ its restriction.

We consider two kinds of classes on $A^*_{\aut(C_0)}$ induced by the sheaf $\omega^{\Gamma}_0$:
\begin{enumerate}
\item the pushforward of polynomials of the Chern class $K:=c_1(\omega^{\Gamma}_0)$;
\item the Chern classes of the pushforward of $\omega_0$.
\end{enumerate}  

The first kind will give us the equivalent of Mumford classes, but on $\M_0^{\leq 3}$ we can only define classes of the second type.

The aim of this section is to compute classes of the first kind for $\Gamma$ with at most three nodes  and to describe them as polynomials in classes of the second kind. If this description is independent from the graph $\Gamma$ then we can define them as elements of $A^*(\M_0^{\leq 3})$.

First of all we define on $\cu^{\Gamma}$ and $\M^{\Gamma}_0$ the following classes
\begin{eqnarray*}
K &:=& c_1(\omega^{\Gamma}_0) \in A^1(\cu^{\Gamma}),\\
\km_i &:=& \Pi_*(K^{i+1}) \in A^{i}(\M^{\Gamma}_0).
\end{eqnarray*}

Such classes $\km_i$ (introduced in \cite{mumenum} for the moduli spaces of stable curves) are called {\it Mumford classes}.

{\bf Mumford classes on $\M_0^0$}

The stack $\M^0_0$ is $\B \Pro Gl_2$.

Let us consider the universal curve
$$
[\Pro ^1 _{\C} / \Pro Gl_2] \xrightarrow{\Pi} \B \Pro Gl_2.
$$
We call $\overline \Pi$ the induced map $\Pro ^1 _{\C} \to \spc$ and $\overline \omega^0_0$  a lifting of $\omega^0_0$ on $\Pro^1_{\C}$.

We notice that $\left ( \overline \omega^0_0 \right ) ^{\vee} = T _ {\overline \Pi}$ ($T_{\Pi}$ is the relative tangent bundle along $\Pi$), we set  $$ K:=c^{\Pro Gl_2}_1(T _{\overline \Pi})= -c^{\Pro Gl_2}_1(\overline \omega^0_0).$$

Furthermore  $\overline \Pi_*(T_{\overline \Pi})=H^0(\Pro^1_{\C}, T_{\overline \Pi})=\sla _2$ seen as adjoint representation of $\Pro Gl_2$.

By applying the equivariant Grothendieck-Riemann-Roch Theorem we obtain
\begin{eqnarray*}
ch(\Pi_*(\omega_0^\vee)) &=& \Pi_*(Td(T_{\Pi})ch(\omega_0^\vee))\\
ch^{\Pro Gl_2}(\sla_2)&=&\Pi_*(Td^{\Pro Gl_2}(T_{\overline \Pi})ch^{\Pro Gl_2}(T_{\overline \Pi}))\\
3 -c^{\Pro Gl_2}_2(\sla _2)&=& \Pi_* \left[ \left( e^{- (K)}\right) \left( \frac{-K}{1- e^{K}}\right) \right]
\end{eqnarray*}
By applying GRR to the trivial linear bundle we obtain:
$$
1= \Pi_* \left[ \frac{-K}{1- e^{K}} \right]
$$
If we subtract the second equation from the first, we obtain
$$
2 -c^{\Pro Gl_2}_2(\sla _2)  = \Pi_* \left[ \left( \frac{1-e^{K}}{e^{K}}\right) \left( \frac{-K}{1- e^{K}}\right) \right] = \Pi_* \left[ -K e^{-K} \right],
$$
from which we get the following
\begin{propos}\label{mumzero}
On $\M^0_0$ we have
\begin{eqnarray*}
\km_0 = -2, \quad \km_2 = 2c^{\Pro Gl_2}_2(\sla _2), \quad \km_1 = \km_3 = 0.
\end{eqnarray*}
\end{propos}

{\bf Mumford classes on strata of singular curves}

Now let us consider the following cartesian diagram (see Section \ref{stratclass})
\begin{equation}\label{diagtildef}
\xymatrix@=3pc{
[\stackrel{\delta + 1}{\coprod} \Pro ^1 _{\C} /  H] \cart \ar@/_3pc/[dd]_{\widetilde F} \ar [r] ^\xi \ar[d]_{\widetilde N} & \widehat{\cu} ^{\Gamma} \ar@/^2pc/[dd]^{F} \ar[d]^{N} \\
[C_0 / H]  \cart \ar [d]_{\widetilde \Pi} \ar[r] &\cu ^{\Gamma} \ar[d]^\Pi\\
\B H \ar[r] _{\phi} & \B  \aut(C_0)
}
\end{equation}
where $N: \widehat{\cu} ^{\Gamma} \to \cu^{\Gamma}$ is the normalization of $\cu^{\Gamma}$.

In the following we call $\widetilde{\omega}^{\Gamma}_0$ the sheaf
$$
\xi^* N^*(\omega^{\Gamma}_0).
$$ 
and $\widetilde{K}:=c_1(\omega_0)$.

\begin{propos} \label{pfkm}
Using the above notation, the Mumford classes $\km_i \in A^i(\B F \ltimes H) \rat$ are described by the following relation
$$
\phi_* \widetilde{F}_*(\widetilde{K}^{i+1}) = \sigma \km_i,
$$
\end{propos}

\begin{dem}
Let us fix an index $i \in \N$. Since the map $N$ is finite and generically of degree 1 we have $N_* N^* = \id$
therefore $\km_i := \Pi_* K^{i+1} = \Pi_*(N_* N^*)K^{i+1} = F_*(N^* K^{i+1})$,
since $F$ is projective and $\phi$ is \'etale we can apply the projection formula and obtain $$\phi^* F_*(N^* K^{i+1})= \widetilde{F}_* \xi ^* (N^* K^{i+1}).$$
Furthermore the map $N \circ \xi$ is finite and so $$\xi^* N^* c_1(\omega^{\Gamma}_0)^{i+1} = c_1(\xi^* N^* (\omega^{\Gamma}_0))^{i+1} = \widetilde K ^{i+1}.$$
We conclude by noting that $\phi_* \phi^*$ is multiplication by $\sigma$.
\end{dem}

\begin{propos}\label{defmumstrat}
Following notation of Proposition \ref{classstratum} and order elements of $\Delta_1$ from $1$ to $\delta_1$, we have
\begin{eqnarray} \label{kmum}
\km_m=-\phi_*\frac{t_1^m+ \dots + t_{\delta_1}^m}{\sigma}.
\end{eqnarray}
\end{propos}

\begin{dem}
From Proposition \ref{pfkm}, we have reduced the problem to computing $\widetilde K$ and then writing pushforward along $\widetilde{F}$. With respect to each component of the stack
$$
\widehat{\cu} ^{\Gamma} = \stackrel{\delta + 1}{\coprod} [\Pro ^1 _{\C} / \aut(C_0)]
$$
we fix coordinates on $\Pro ^1 _{\C}$ which are compatible with coordinates chosen on $C$. By forgetting the action of $\aut(\Gamma)$ we keep the same system of coordinates on $\stackrel{\delta + 1}{\coprod} \Pro ^1 _{\C}$ when we consider $\stackrel{\delta + 1}{\coprod} [\Pro ^1 _{\C} / H]$.

Giving a bundle on a quotient stack $[X/G]$ is equivalent to giving a bundle  $U \to X$ equivariant for the action of $G$. 
By abuse of notation we still call $\widetilde{\omega}^{\Gamma}_0$ any lifting of $\widetilde{\omega}^{\Gamma}_0$ on $\stackrel{\delta + 1}{\coprod} \Pro ^1 _{\C}$.

In order to make computations we need to render explicit the action of $H:=E^{\Delta_1} \times \Gm^{\Delta_2}$ on
$$
\widetilde{F}_*(\widetilde{\omega}_0)= H^0 \left( \stackrel{\delta + 1}{\coprod} \Pro ^1 _{\C}, \widetilde{\omega}^{\Gamma}_0 \right ).
$$
 
Set $\Delta:=\Delta_1 \cup \Delta_2$, from the inclusion $\Gm \to E$ we have a cartesian diagram of stacks
\begin{equation*}
\xymatrix{
\B (\Gm)^{\Delta} \cart \ar[r]^-{\phi} \ar[d]^{\Psi} &\B (\aut(\Gamma) \ltimes (\Gm)^{\Delta}) \ar[d]^{\Psi}\\
\B H \ar[r]^{\phi} &\B \aut (C_0)
}
\end{equation*}
Roughly speaking we can say that we obtain the stacks in the top row by fixing the point $0$ on components with a node. The functor $\Psi$ forgets these points.
We have the following ring isomorphisms 
\begin{eqnarray*}
A^*_{H} &\xrightarrow{\Psi^*}& A^*_{(\Gm)^{\Delta}} \\
 A^*_{\aut(C_0)} &\xrightarrow{\Psi^*}& A^*_{\aut(\Gamma) \ltimes (\Gm)^{\Delta}}.
\end{eqnarray*}
We have defined the classes $t_{\Delta_1}, r_{\Delta_2}$ in $A^*_H$. By using the same notation of Section \ref{stratclass}, the map $\Psi^*$ is the identity on $r_{\Delta_2}$. For each vertex $P$ in $\Gamma$ such that $e(P)=1$ the map $\Psi^*$ sends $t_P$ to $c^{\Gm}_1(H^0(\Pro^1_{\C},\Of(z_\infty)))$ of the same component, which we still call $t_P$.
Consequently, with reference to the following diagram
\begin{equation*}
\xymatrix@=3pc{
[\stackrel{\Delta}{\coprod} \Pro ^1 _{\C} /  (\Gm)^{\Delta}] \ar[d]_{\widetilde{F}} \\
\B (\Gm)^{\Delta} \ar[r] _-{\phi} & \B  (\aut(\Gamma) \ltimes (\Gm)^{\Delta})
}
\end{equation*}
we have reduced the problem to consider the action of $(\Gm)^{\Delta}$ on
$$
\widetilde{F}_*(\widetilde{\omega}^{\Gamma}_0)= H^0 \left( \stackrel{\Delta}{\coprod} \Pro ^1 _{\C}, \widetilde{\omega}^{\Gamma}_0 \right )
$$
where, again with abuse of notation, we call $\widetilde{\omega}^{\Gamma}_0$ the sheaf $\Psi^*(\widetilde{\omega}^{\Gamma}_0)$.

The map $\widetilde{F}$ is the union of maps	
$$
\widetilde{F}_P:= [\Pro^1 _{\C}/ (\Gm)^\Delta] \to \B (\Gm)^\Delta
$$
where only the component of $(\Gm)^{\Delta}$ corresponding to $P \in \Delta$ does not acts trivially on $\Pro^1_{\C}$ and the action is
\begin{eqnarray*}
\az: \Gm \times \Pro^1_{\C} &\to& \Pro^1_{\C}\\
(\lambda, [X_0,X_1]) &\mapsto& [X_0, \lambda X_1].
\end{eqnarray*}
Set
$$
P_0:=[0,1], \quad P_1:=[1,1], \quad P_{\infty}:= [1,0].
$$
We have that $\mathcal F := (\widetilde{\omega}_0)^\vee$ is the sheaf on  $\stackrel{\Delta}{\coprod}\Pro ^1 _{\C}$ such that restricted:
\begin{itemize}
\item to the components with one node it is $\mathcal F_P:=(\omega \otimes \Of (z_{\infty}))^\vee$
\item to the components with two nodes it is $\mathcal F_P:=(\omega \otimes \Of (z_{\infty} + z_0))^\vee$
\item to the components with three nodes it is $\mathcal F_P:=(\omega \otimes \Of (z_{\infty} + z_0 + z_1))^\vee$
\end{itemize}
where $\omega$ is the canonical bundle on $\Pro^1_{\C}$.

In the following, Chern classes will be equivariant for the action of $\Gm$. For each $P \in \Delta$ let us indicate on each component
\begin{eqnarray*}
K_P &:=& c^{\Gm}_1(\omega) \in A^1_{\Gm}(\Pro^1_{\C})\\
R_P &:=& c^{\Gm}_1(\Of(z_{\infty})) \in A^1_{\Gm}(\Pro^1_{\C})\\
Q_P &:=& c^{\Gm}_1(\Of(P_0)) \in A^1_{\Gm}(\Pro^1_{\C})
\end{eqnarray*}

We have
\begin{eqnarray*}
c^{\Gm}_1(\F_{\Delta_1})&=& - {K_{\Delta_1}} - R_{\Delta_1}\\
c^{\Gm}_1(\F_{\Delta_2})&=& - {K_{\Delta_2}} - R_{\Delta_2} - Q_{\Delta_2}.
\end{eqnarray*}
In order to determine Mumford classes it is then necessary to compute the pushforward classes
\begin{eqnarray*}
&& \widetilde{F}_{\Delta_1*} (- {K_{\Delta_1}} - R_{\Delta_1})^h \text{ and}\\
&& \widetilde{F}_{\Delta_2*} (- {K_{\Delta_2}} - R_{\Delta_2} - Q_{\Delta_2})^h
\end{eqnarray*}
for every natural $h$.

Let us start with computing the push-forward along $\widetilde F _{\Delta}$ of every power of $K_{\Delta}=c^{\Gm}_1(\omega )$.

We have
$$
\widetilde F_{\Delta_1*}(\omega^\vee)=H^0(\Pro ^1 _\C, \omega^\vee).
$$
Now it is necessary to determine the action of $(\Gm)$
on global sections of  $\omega ^\vee$. Let $z:= X_1 / X_0$ be the local coordinate around $z_{\infty}$, the global sections of $\omega^{\vee}$ are generated as a vectorial space by $\hol, z \hol, z^2 \hol$. Relatively to this basis, the action of $\Gm$ is given by
\begin{eqnarray*}
\az: \Gm \times H^0(\Pro^1_{\C},\omega^\vee) &\to& H^0(\Pro^1_{\C}, \omega^\vee)\\
(\lambda, a\hol + b z \hol + c z^2 \hol) &\mapsto& (a \lambda \hol + b z \hol + c \frac{1}{\lambda} z^2 \hol)
\end{eqnarray*}
the multiplicity of the action is therefore $(1,0,-1)$, so the Chern character of $\widetilde F _{\Delta*}(\omega^\vee)$ is $e^{t_{\Delta}} + 1 + e^{-t_{\Delta}}$.

Let us observe that $\omega^{\vee}$ is the tangent bundle relative to $\widetilde F _{\Delta}$, so by applying the GRR Theorem we get
$$
1 + e^{t_{\Delta}} + e^{-t_{\Delta}} = \widetilde F_{\Delta*} \left[ \left( e^{- (K_{\Delta})}\right) \left( \frac{-K_\Delta}{1- e^{K_\Delta}}\right) \right]
$$
By applying GRR to the trivial linear bundle we obtain:
$$
1= \widetilde F_{\Delta*} \left[ \frac{-K_\Delta}{1- e^{K_\Delta}} \right]
$$
If we subtract the second equation from the first, we obtain
$$
e^{t_\Delta} + e^{-t_\Delta}  = \widetilde F_{\Delta*} \left[ \left( \frac{1-e^{K_\Delta}}{e^{K_\Delta}}\right) \left( \frac{-K_\Delta}{1- e^{K_\Delta}}\right) \right] = \widetilde F_{\Delta*} \left[ -K_\Delta e^{-K_\Delta} \right],
$$
from this, by distinguishing between even and odd cases, it follows that

\begin{eqnarray*}
\widetilde F_{\Delta*}K_\Delta^{2h}&=&0\\
\widetilde F_{\Delta*}K_\Delta^{2h+1}&=& -2t_\Delta^{2h}
\end{eqnarray*}

Let us notice that there exist two equivariant sections of $\widetilde F_{\Delta}$ given by the fixed points $z_0$ and $z_\infty$, that we will call respectively $s_0$ and $s_\infty$
\begin{equation*}
\xymatrix{
\B (\Gm)^{\Delta} \ar@<0.5ex>[r]^-{s_0} \ar@<-0.5ex>[r]_-{s_{\infty}} & [\Pro^1_{\C}/ (\Gm)^{\Delta} ]
}
\end{equation*}
From the self intersection formula we have 
\begin{eqnarray*}
&&s_0^*(Q_\Delta)=-t_\Delta,  \; s_0^*(R_\Delta)=0, \; s_0^*(K_\Delta)=t_\Delta\\
&&s_\infty^*(Q_\Delta)=0, \; s_\infty^*(R_\Delta)=t_\Delta, \; s_\infty^*(K_\Delta)=-t_\Delta;
\end{eqnarray*}
then by applying the projection formula (see \cite{ful} p.34)
for every cycle $D \in A^*_{(\Gm)^\Delta}(\Pro^1)$ we have
\begin{eqnarray*}
&& Q_\Delta \cdot D = s_{0*}(1) \cdot D = s_{0*}(s_0^*D)\\
&& R_\Delta \cdot D = s_{\infty*}(1) \cdot D = s_{\infty*}(s_\infty^*D).
\end{eqnarray*} 
With reference to the previous relations, the Mumford classes are determined on every component (we separate between even and odd cases)
\begin{eqnarray*}
\widetilde F_{\Delta *}(-K_\Delta -R_\Delta)^{2h} &=& \widetilde F^\Delta_*(K_\Delta)^{2h} + 
					\sum_{a=1}^{2h} \binom{2h}{a}				
					\widetilde F_{\Delta*}(R_\Delta^{a}\cdot K_\Delta^{2h-a})\\
				&=& 	\sum_{a=1}^{2h} \binom{2h}{a}	
					(\widetilde F_{\Delta*}s_{\infty*})s_\infty^*(R_\Delta^{a-1}\cdot K_\Delta^{2h-a})\\
				&=& 		\sum_{a=1}^{2h} \binom{2h}{a}	
						(-1)^{a-1}t_\Delta^{2h-1} = -t_\Delta^{2h-1}
\end{eqnarray*}
and in a completely analogous way, we have the following relations
\begin{eqnarray*}
&&\widetilde F_{\Delta_1*}(-K_{\Delta_1} -R_{\Delta_1})^{2h + 1} = t_i^{2h}\\
&&\widetilde F_{\Delta_2*}(-K_{\Delta_2} -R_{\Delta_2} -Q_{\Delta_2})^{h}=0.
\end{eqnarray*}
The last relation allows us to ignore components with two nodes.
By following the notation of Proposition \ref{pfkm} we notice that
$$
\widetilde F _* \widetilde K ^ {m+1} = \sum_{P \in \Delta_1} \widetilde F_{P*}(-1)^{m+1}(-K_P -R_P)^{m+1} 
$$
and consequently, if we order elements of $\Delta_1$ from $1$ to $\delta_1$, we have

\begin{eqnarray}
\km_m=-\phi_*\frac{t_1^m+ \dots + t_{\delta_1}^m}{\sigma}.
\end{eqnarray}

\end{dem}

{\bf Definition of Mumford classes on $\M_0^{\leq 3}$}
Now we consider the universal curve
$$
\cu^{\leq 3} \xrightarrow{\Pi} \M_0^{\leq 3}.
$$
We have seen in \cite{fulg2} that the pushforward $\Pi_* \left ( \omega^{\leq 3}_0 \right )^{\vee}$ is a well defined rank three vector bundle. Consequently from \cite{kre} Section 3.6 we have that $c_i(\Pi_* \left ( \omega^{\leq 3}_0 \right )^{\vee})=0$ for $i>3$.
We fix the following notation
\begin{eqnarray*}
 \cm_1:=  c_1(\Pi_*\left ( \omega^{\leq 3}_0 \right )^{\vee}), \quad \cm_2:= c_2(\Pi_* \left ( \omega^{\leq 3}_0 \right )^{\vee}), \quad \cm_3:= c_3(\Pi_* \left ( \omega^{\leq 3}_0 \right )^{\vee}).
\end{eqnarray*}
We still call $\cm_1, \cm_2, \cm_3$ their restriction to each stratum of $\M^{\leq 3}_0$.

\begin{defi}\label{defmumpol}
We define in $A^*(\M^{\leq 3}_0)$ Mumford classes $\km_1, \km_2, \km_3$ as follows
\begin{eqnarray*}
\km_1 :=  -\cm_1, \quad \km_2  :=  2 \cm_2 - \cm^2_1, \quad \km_3 := -\cm^3_1 + 3 \cm_1 \cm_2 -3\cm_3.
\end{eqnarray*}
\end{defi}

\begin{oss}
From Proposition \ref{injleqtre}, the Chow ring $A^*(\M^{\leq 3}_0) \rat$ injects into the product of $A^*(\M^{\Gamma}_0)\rat$ over trees with at most three edges. Consequently in order to verify that the above is a good definition we only need to prove that the restrictions of Mumford classes to each stratum are the given polynomials.
\end{oss}

\begin{propos}\label{finmum}
Let $\Gamma$ be a tree with at most three edges. Set 
\begin{eqnarray*}
 \cm_1(\Gamma) &:=& c_1(\Pi_*\left ( \omega^{\Gamma}_0 \right )^{\vee})\\
 \cm_2(\Gamma) &:=& c_2(\Pi_* \left ( \omega^{\Gamma}_0 \right )^{\vee})\\
 \cm_3(\Gamma) &:=& c_3(\Pi_* \left ( \omega^{\Gamma}_0 \right )^{\vee}).
\end{eqnarray*}
and
\begin{eqnarray*}
&& \nm_1(\Gamma) :=  \cm_1(\Gamma), \; \nm_2(\Gamma)  :=  \cm^2_1(\Gamma) -2 \cm_2(\Gamma), \\
&& \nm_3(\Gamma) := \cm^3_1(\Gamma) - 3 \cm_1(\Gamma) \cm_2(\Gamma) +3\cm_3(\Gamma).
\end{eqnarray*}

Then we have
\begin{eqnarray*}
\km_1 =  -\nm_1(\Gamma) \quad \km_2  =  -\nm_2(\Gamma) \quad \km_3 =  -\nm_3(\Gamma)
\end{eqnarray*}
\end{propos}

\begin{dem}
We consider first of all the case $\M^{\Gamma}_0= \M^0_0$. Here we have
\begin{eqnarray*}
&&\nm_1(\Gamma) = c^{\Pro Gl_2}_1(\sla_2)=0\\
&&\nm_2(\Gamma) =-2c^{\Pro Gl_2}_2(\sla_2)\\
&&\nm_3(\Gamma)=0
\end{eqnarray*}
and we can conclude because on $\M^0_0$ we have that $\km_1=\km_3=0$ and $\km_2=2c^{\Pro Gl_2}_2(\sla_2)$ (Proposition \ref{mumzero}).

Now let us consider any other tree $\Gamma$ with at most three edges.  With reference to diagram (\ref{diagtildef}),
as Chern classes commutes with base changing, we have
$$
\phi^*c_i(\Pi_*\left( \omega^{\Gamma}_0 \right)^{\vee}) = c_i(\widetilde \Pi_* \left( \widetilde \omega^{\Gamma}_0 \right) ^{\vee}).
$$
Given $\Gamma$ we order elements of $\Delta_1$ from $1$ to $\delta_1$. By putting together the above relation and the equation (\ref{kmum}), we reduce to show
$$
ch(\widetilde \Pi_* \left( \widetilde \omega^{\Gamma}_0 \right)^{\vee})= 3 + \sum^{\infty}_{m=1} \frac{t^m_1 + \dots + t^m_{\delta_1}}{m!}.
$$
We recall that the universal curve $\widetilde \cu^{\Gamma}$ on $\widetilde \M_0^{\Gamma}= \B (E^{\Delta_1} \times \Gm^{\Delta_2})$ is the quotient stack $[C_0/ E^{\Delta_1} \times \Gm^{\Delta_2}]$.

Since we have a morphism $\Psi: \B (\Gm)^{\Delta} \to \widetilde \M_0^{\Gamma}$
such that $\Psi^*: A^*(\widetilde \M_0^{\Gamma}) \to A^*_{(\Gm)^{\Delta}}$ is an isomorphism, with reference to the cartesian diagram
\begin{equation*}
\xymatrix{
[C_0 / \Gm^{\Delta}] \cart \ar[r]^-\Psi \ar[d]_{\widetilde \Pi} & [C_0/E^{\Delta_1} \times \Gm^{\Delta_2}] \ar[d]^{\widehat \Pi}\\
\B \Gm^{\Delta} \ar[r]^-\Psi & \B E^{\Delta_1} \times \Gm^{\Delta_2}
}
\end{equation*}
we reduce to consider the pullback sheaf $\Psi^* \left( \widetilde \omega^{\Gamma}_0 \right)^\vee$ which we still call  $\left( \widetilde \omega^{\Gamma}_0 \right)^\vee$ as its lifting to $C_0$.
Standard arguments involving Serre duality, show that $H^1(C_0, \left( \widetilde \omega^{\Gamma}_0 \right)^\vee)=0$, so we have
\begin{equation*}
ch(\widetilde \Pi_* \left( \widetilde \omega^{\Gamma}_0 \right)^{\vee})= ch^{\Gm^\Delta}(H^0(C_0, \left( \widetilde \omega^{\Gamma}_0 \right)^\vee))
\end{equation*}
On curves of topological type
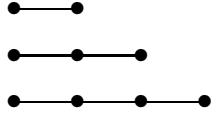
\begin{eqnarray*}
&& \vcenter{\unoa}\\
&& \vcenter{\duea}\\
&& \vcenter{\trea}
\end{eqnarray*}
we have that global sections of $ \left( \widetilde \omega^{\Gamma}_0 \right)^{\vee}$ are sections of $\Of_{\Pro^1_{\C}}(1)$ on extremal components and $\Of_{\Pro^1_{\C}}$ on the other components, which agree on nodes. Since $H^0(\Pro^1_\C,\Of_{\Pro^1_{\C}}(0))= \C$,  global sections of $ \left( \widetilde \omega^{\Gamma}_0 \right)^{\vee}$ are sections of $\Of_{\Pro^1_{\C}}(1)$ on the two extremal components which are equal on nodes. On each extremal component we fix coordinates $[X'_0,X'_1]$ and $[X''_0, X''_1]$. Sections on $\Of_{\Pro^1_{\C}}(1)$ are linear forms
\begin{eqnarray*}
&& a_1 X'_0 + b_1 X'_1\\
&& a_2 X''_0 + b_2 X''_1
\end{eqnarray*}
which agree at $z_{\infty}=[1,0]$. This happens if and only if  $a_1=a_2$. We have only $(\Gm)^{\Delta_1}$ which does not acts trivially on $H^0(C_0, \left( \widetilde \omega^{\Gamma}_0 \right)^{\vee})$ and the action is
\begin{eqnarray*}
\az: (\Gm \times \Gm) \times H^0(C_0, \left( \widetilde \omega^{\Gamma}_0 \right)^{\vee}) &\to& H^0(C_0, \left( \widetilde \omega^{\Gamma}_0 \right)^{\vee})\\
(\lambda_1,\lambda_2), (a, b_1, b_2 ) &\mapsto& ( a, \lambda_1 b_1, \lambda_2 b_2)
\end{eqnarray*}
so the Chern character of $\widetilde \Pi_* \omega^{\vee}_0$ is
$$
1 + e^{t_1} + e^{t_2} = 3 + \sum^{\infty}_{m=1} \frac{t^m_1+ t^m_2}{m!}
$$
as we claimed.

The last case to consider is when $\Gamma$ equals to
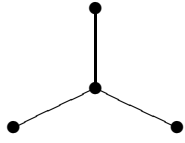
\begin{eqnarray*}
&& \vcenter{\treb}\\
\end{eqnarray*}
We have that global sections of $\left( \widetilde \omega^{\Gamma}_0 \right)^{\vee}$ are sections of $\Of_{\Pro^1_{\C}}(1)$ on extremal components and $\Of_{\Pro^1_{\C}}(-1)$ on the central component which agree on nodes. Since $H^0(\Pro^1_\C,\Of_{\Pro^1_{\C}}(-1))= 0$,  global sections of $\left( \widetilde \omega^{\Gamma}_0 \right)^{\vee}$ are sections of $\Of_{\Pro^1_{\C}}(1)$ on the three extremal components which are zero on nodes. On each extremal component we fix coordinates $[X'_0,X'_1]$, $[X''_0,X''_1]$ and $[X'''_0, X'''_1]$. Sections on $\Of_{\Pro^1_{\C}}(1)$ are linear forms
$$ 
a_1 X'_0 + b_1 X'_1 \quad a_2 X''_0 + b_2 X''_1 \quad a_3 X'''_0 + b_3 X'''_1
$$
which are zero on nodes if and only if $a_1=a_2=a_3=0$. We have that only $(\Gm)^{\Delta_1}$ does not acts trivially on $H^0(C_0,\left( \widetilde \omega^{\Gamma}_0 \right)^{\vee})$ and the action is
\begin{eqnarray*}
\az: (\Gm \times \Gm \times \Gm) \times H^0(C_0,\left( \widetilde \omega^{\Gamma}_0 \right)^{\vee}) &\to& H^0(C_0,\left( \widetilde \omega^{\Gamma}_0 \right)^{\vee})\\
(\lambda_1,\lambda_2, \lambda_3), (b_1, b_2, b_3 ) &\mapsto& ( \lambda_1 b_1, \lambda_2 b_2, \lambda_3 b_3)
\end{eqnarray*}
so the Chern character of $\widetilde \Pi_* \omega^{\vee}_0$ is
$$
e^{t_1} + e^{t_2} + e^{t_3} = 3 + \sum^{\infty}_{m=1} \frac{t^m_1+ t^m_2 + t^m_3}{m!}
$$
as we claimed.
\end{dem}

\end{document}